\documentclass[12pt]{article}
\usepackage{amsmath,amssymb,amsthm,latexsym,graphicx}
\newtheorem{theorem}{Theorem}
\newtheorem{lemma}[theorem]{Lemma}
\newtheorem{remark}[theorem]{Remark}
\newtheorem{corollary}[theorem]{Corollary}
\newtheorem{definition}[theorem]{Definition}
\newcommand{\R}{\mathbb{R}}

\newcommand{\bd}{\begin{definition}}
\newcommand{\ed}{\end{definition}}
\newcommand{\bt}{\begin{theorem}}
\newcommand{\et}{\end{theorem}}
\newcommand{\bl}{\begin{lemma}}
\newcommand{\el}{\end{lemma}}
\newcommand{\bc}{\begin{corollary}}
\newcommand{\ec}{\end{corollary}}
\newcommand{\bcon}{\begin{conjecture}}
\newcommand{\econ}{\end{conjecture}}
\newcommand{\br}{\begin{remark}}
\newcommand{\er}{\end{remark}}
\newcommand{\bp}{\begin{proposition}}
\newcommand{\ep}{\end{proposition}}
\newcommand{\be}{\begin{equation}}
\newcommand{\ee}{\end{equation}}
\newcommand{\bef}{\begin{figure}}
\newcommand{\eef}{\end{figure}}
\newcommand{\bea}{\begin{eqarray}}
\newcommand{\eea}{\end{eqarray}}
\newcommand{\ba}{\begin{array}}
\newcommand{\ea}{\end{array}}
\begin{document}
\title{On a problem concerning affine-invariant points of convex sets}
\author{Peter Kuchment\\Texas A\&M University\\kuchment@math.tamu.edu}
\date{}
\maketitle
\begin{abstract}
This text is a somewhat reformated\footnote{ E.g., some statements that were not as such in the original paper, are given the names ``Corollary'' or ``Theorem.''} translation of the old and practically inaccessible paper: P. Kuchment,  On the question of the affine-invariant points of convex bodies,
(in Russian),  Optimizacija  No. 8(25)  (1972), 48--51, 127. MR0350621. There partial solutions of some old problems of B.~Gr\"unbaum concerning affine-invariant points of convex bodies were obtained. The main restriction, due to which the solution was incomplete, was the compactness restriction on the group of linear transformations involved. It was noticed recently by O.~Mordhorst (arXiv:1601.07850) that a simple additional argument allows one to restrict the consideration to the bodies whose John's ellipsoid is a ball and consequently to a compact group case. This in turn extends the result to the complete solution of B.~Gr\"unbaum's problem (see also some partial progress in the previous recent works by M.~Meyer, C.~Sch\"utt, and E.~Werner \cite{MSW,MSW2}).

The article translated below contains three main statements: existence of invariant points in the case when the linear parts of affine transformations are taken from a compact group, the resulting from this proof of B.~Gr\"unbaum's conjecture for the case of similarity invariant points, and an example precluding the set of all such points from being the convex hull of finitely many of them.
\end{abstract}
Let $\R^n$ be the standard $n$-dimensional real vector space. We denote by $W_n$ the set of all compact convex subsets in $\R^n$ with non-empty interiors (we will call them \emph{convex bodies}).The set $W_n$ is equipped with the Hausdorff metric topology \cite{Grun}. Let $G$ be a compact subgroup of the group $GL(n)$ of invertible linear transformations of $\R^n$. We denote by $\Pi$ the group of affine transformations generated by shifts and elements of $G$.

\bd
A \textbf{$\Pi$-invariant point} is a continuous mapping $p:W_n \mapsto \R^n$ such that
\be
p(A(K))=A(p(K))
\ee
for any $K\in W_n$ and $A\in\Pi$.
\ed

\bd
A $\Pi$-invariant point $p$ such that $p(K)\in K$ for all $K\in W_n$, is called a \textbf{$\Pi$-invariant point of the body.}
\ed
In what follows, we will use the abbreviations \textbf{$\Pi$-i point} and \textbf{$\Pi$-i point of the body} correspondingly.

Let $K\in W_n$. We denote by $\Pi_K$ the subgroup of $\Pi$ consisting of all automorphisms of $K$ (i.e., those leaving $K$ invariant).
\bt\indent
\begin{enumerate}
\item Let $x_0$ be a $\Pi_{K_0}$-invariant point, located in the interior $Int\, K_0$ of a convex body $K_0$. Then there exists a $\Pi$-invariant point of the body $p:W_n\mapsto\R^n$ such that $p(K_0)=x_0$.
\item Let $x_0$ be an arbitrary  $\Pi_{K_0}$-invariant point. Then there exists a $\Pi$-invariant point $p:W_n\mapsto\R^n$ such that $p(K_0)=x_0$.

        In other words, given a convex body $K$, the set of all $\Pi_K$-invariant points in $\R^n$ coincides with the set of values $p(K)$ of all $\Pi$-invariant points $p$.
\end{enumerate}
\et
\begin{proof} We will restrict ourselves to the proof of the first statement, since the second one is proven analogously. We will need to describe some new objects. All required properties of equivariant bundles can be found in \cite{Husem}.

We define the set $W_n^0$ of classes of equivalence $k=\{K\}$ in $W_n$, with respect to the shift equivalence. A metric on $W_n^0$ can be defined as follows:
\be
\rho(k_1,k_2)=\inf\limits_{K_j\in k_j}\rho(K_1,K_2).
\ee
Here $\rho$ on the right denotes the Hausdorff distance between convex bodies.
The general linear group $GL(n)$, and thus also its subgroup $G$, act naturally on $W_n^0$:
\be
g(\{K\})=\{g(K)\} \mbox{ for } g\in GL(n).
\ee
Let us consider the set $\mathcal{E}:=W_n^0\times \R^n$ as a trivial fiber-bundle over $W_n$, with respect to the natural projection onto the first factor. One can introduce on it the structure of a $G$-bundle:
\be
g(k,x):=(g(k), g(x)).
\ee
We can now define an equivariant sub-bundle $\mathcal{K}\subset\mathcal{E}$ as follows: let $K\in k \in W_n^0$, then the fiber of $\mathcal{K}$ over $k$ is the body $K\subset \R^n$, after shifting it so that its center of mass coincides with the origin $0$. Due to the continuous dependence of the center of mass on the convex body \cite{Grun}, $\mathcal{K}$ is indeed a bundle.

In these new terms, a $\Pi$-invariant point (a $\Pi$-invariant point of the body)is a continuous equivariant section of the bundle $\mathcal{E}$ (respectively, of $\mathcal{K}$). Now the first statement of the theorem can be rephrased as follows:
\begin{center}
 \emph{ Let $x_0\in Int \,K_0$ be $G_{K_0}$-invariant (the stationary subgroup $G_K$ is defined analogously to $\Pi_K$). Then it can be extended to a continuous equivariant section $p$ of $\mathcal{K}$ such that $p(\{K_0\})=x_0$.}
\end{center}
This is the form in which we will be proving it. It is known \cite{Grun} that $W_n^0$ is locally compact. We choose a compact neighborhood $\mathcal{V}$ of $\{K_0\}$. Due to compactness of $G$, the orbit $G(\mathcal{V})$ is also compact. Let us restrict the bundle $\mathcal{E}$ to $G(\mathcal{V})$. Over the single orbit $G(\{K\})$ we can define an equivariant section\footnote{\textbf{Remark during translation}: $G_K$-invariance of $x_0$ is used.} $\delta$ as $\delta(g(\{K\}))=g(x_0)$. It is known \cite{Husem} that, due to the closedness of $G(\{K\})$ and compactness of $G(\mathcal{V})$, $\delta$ can be extended to a continuous equivariant section (which we denote with the same letter $\delta$) over $G(\mathcal{V})$.

Since $x_0$ belongs to the interior of $K_0$, there is a $G$-invariant neighborhood $\mathcal{U}$ of the orbit, where $\delta(F)\in Int\, F$ for any $F\in\mathcal{U}$. Consider covering of $W_n^0$ with two $G$-invariant open sets:
$$
U_1:=\mathcal{U},\,\, U_2:=W^0_n\setminus G(\{K_0\}).
$$
Due to metrizability of $W_n^0$, there exists a subordinate partition of unity $\phi_1(k), \phi_2(k)$. We can assume these two functions being $G$-invariant. Otherwise, we can redefine them as follows:
\be
 \tilde{\phi}_i(\{K\}) = \int_G \phi_i(g(\{K\}))d\mu(g),
\ee
where $\mu(g)$ is the normalized Haar measure\footnote{\textbf{Remark during translation:} Here compactness of $G$ is used again.} on $G$.

We can now define
\be
p(\{K\}):= \phi_1(\{K\})\delta(\{K\})+\phi_2(\{K\})\theta(\{K\}),
\ee
where $\theta(\{K\})$ is the zero section of $\mathcal{K}$, i.e. the center of mass of $K$.

It is clear that $p$ is a continuous equivariant section of $\mathcal{K}$ and $p(\{K_0\})=x_0$. This proves the first statement of the theorem. The second one, as it has been mentioned before, is proven by the same method, but even easier (no restriction $p(K)\in K$ to obey).
\end{proof}

In \cite{Grun}, the questions were asked both for affine- and \emph{similarity- invariant} points. The first case corresponds to $G=GL(n)$, and the second one to the group $G$ of linear transformations preserving orthogonality. Both these groups are non-compact and thus are not covered by our theorem. However, in the second case, of similarity invariant points, one can reduce consideration to a compact group case considered above. Indeed, one can assume that the volume of $K_0$ is equal to $1$. We can construct the required section $p$ for the class of bodies of unit volume. Indeed, in this case we deal with the compact orthogonal group $O(n)$ and the theorem applies. Then, for a body $F$ of an arbitrary volume $V$, we can define the required section as follows:
\be
p(\{F\}):=V^{1/n}p(\{V^{-1/n}F\}).
\ee
\bc
The Gr\"unbaum's problem has positive solution for similarity invariant points.
\ec
However, for the $GL(n)$ case this reduction to the unit volume bodies does not help, since one arrives to the non-compact group of matrices of unit determinant. The difficulty in the $GL(n)$ case arises when trying to construct an equivariant section of $\mathcal{E}|_{G(\mathcal{V})}$, as it was done above. However, the statement of the theorem probably still holds.

Another question asked in \cite{Grun}:
\begin{center}
\emph{is the set of all similarity-invariant points of the body (clearly a convex set) the convex hull of finitely many such points?}
\end{center}

The answer to this question happens to be negative:

\bt
The set of all similarity-invariant points of the body is not a convex hull of finitely many such points.
\et
\begin{proof} Indeed, let us construct a counterexample. Consider the coordinate plane $x_1=0$ in $\R^n$ and a convex $(n-1)$-dimensional body $K^\prime$ in it, such that its only similarity automorphism is the identity. Let us now construct an $n$-dimensional body $K$ in such a way, that it consists of
\begin{enumerate}
  \item points $(x_1,y)$, where $0<x_1< 1$ and $y(1-x_1)\in K^\prime$,
  \item points $(x_1,y)$, where $0>x_1> -1$ and $y(1+x_1)\in K^\prime$,
  \item points $(1,0)$ and $-1,0$.
\end{enumerate}
This is what is called in topology ``suspension'' over $K^\prime$.

It is clear now that the set of similarity-invariant points of $K$ contains $Int\, K^\prime$ and is contained in $K^\prime$. As soon as $K^\prime$ is not a polyhedron, this is a counterexample.\end{proof}

The author expresses his gratitude to his adviser S. G. Krein, as well as to M. G. Zaidenberg and A. A. Pankov for discussions and comments.

\end{document}